\documentclass[11pt]{article}

\usepackage{amsmath}   
\usepackage{amsthm}          
\usepackage{relsize}

\theoremstyle{plain}

\theoremstyle{definition}

\makeatletter
\@addtoreset{equation}{chapter}

\makeatother

\def\eqref#1{(\ref{#1})}
\def\dsp{\displaystyle}
\def\Frac#1#2{\frac
{
 {\raise.6ex
 \hbox{$\displaystyle#1$}}
}
{
 {\lower.6ex
 \hbox{$\displaystyle#2$}}
 }
}

\numberwithin{equation}{section}

\def\bigOxe{\sqcup \kern-2.3mm \sqcap}

\makeatother

\def\dsp{\displaystyle}
\def\Frac#1#2{\frac
{
 {\raise.6ex
 \hbox{$\displaystyle#1$}}
}
{
 {\lower.6ex
 \hbox{$\displaystyle#2$}}
 }
}

\def\sign{{\rm sign}}

\input epsf

\def\CHFs#1#2#3{
{}_1F_1\left({a};{c};{z}\right)
}

\def\erfc{{\rm erfc}}

\def\tfrac#1#2{{{\lower.6ex
\hbox{$\scriptstyle#1$}}\over 
{\raise.7ex
\hbox{$\scriptstyle#2$}}}}

\def\sign{{\rm sign}}

\def\erfc{{{\rm erfc}}}

\def\erfc{{\rm erfc}}

\def\tfrac#1#2{{{\lower.6ex
\hbox{$\scriptstyle#1$}}\over 
{\raise.7ex
\hbox{$\scriptstyle#2$}}}}

\def\insil#1{}

\begin{document}
 \title{Efficient algorithms for the inversion of the cumulative central beta distribution}

\author{
A. Gil\\
Departamento de Matem\'atica Aplicada y CC. de la Computaci\'on.\\
ETSI Caminos. Universidad de Cantabria. 39005-Santander, Spain.\\ 
\and
J. Segura\\
        Departamento de Matem\'aticas, Estad\'{\i}stica y 
        Computaci\'on,\\
        Univ. de Cantabria, 39005 Santander, Spain.\\
\and
N.M. Temme\\
  IAA, 1391 VD 18, Abcoude, The Netherlands\footnote{Former address: Centrum Wiskunde \& Informatica (CWI), 
        Science Park 123, 1098 XG Amsterdam,  The Netherlands}\\
}

\date{\ }

\maketitle
\begin{abstract}
Accurate and efficient algorithms for the inversion of the cumulative central beta distribution  are
described. The algorithms are based on the combination of a fourth-order
fixed point method with good non-local convergence properties (the Schwarzian-Newton method), asymptotic inversion methods
and sharp bounds in the tails of the distribution function. 
\end{abstract}

\section{Introduction}\label{sec:int}

The cumulative central beta distribution (also known as the incomplete beta function) is defined by 

\begin{equation}
\label{int:01}
I_x(p,q)=\Frac{1}{B(p,q)}\displaystyle\int_0^x t^{p-1}(1-t)^{q-1}\,dt,
\end{equation}
where we assume that $p$ and $q$ are real positive parameters and $0\le x \le 1$.
$B(p,q)$ is the Beta function

\begin{equation}
\label{int:02}
B(p,q)=\Frac{\Gamma(p)\Gamma(q)}{\Gamma(p+q)}.
\end{equation}

From the integral representation in (\ref{int:01}) it is easy to check the following
relation:

\begin{equation}
\label{eq:int03}
I_x(p,q)=1-I_{1-x}(q,p).
\end{equation}

 In this paper we describe algorithms for solving the equation

\begin{equation}
\label{eq:int04}
I_x(p,q)=\alpha,\,\,0<\alpha<1,
\end{equation}
with $p$, $q$ given fixed real positive values. In statistical terms, we are computing the quantile function
for $I_x(p,q)$. 

The beta distribution ia a standard and widely used statistical
distribution which has as particular cases other important distributions like the Student's distribution,
the F-distribution and the binomial distribution. Therefore, the computational schemes for inverting the central beta distribution can be used
 to compute percentiles for other distributions related to the beta. For an example for the F-distribution see \cite{Aber:1993:PFP}.

The quantile function is useful, for instance, for the generation of random variables following the beta distribution density. In some
Monte Carlo simulations the generation of such random variables are required and a massive number of inversions of the beta cumulative distribution
are needed. Therefore, it is important to construct methods as reliable and efficient as possible.

Existing algorithms use some simple initial approximations which are improved by iterating with the Newton method.
In particular, this is the approach used in the inversion method of the statistical software package {\bf R}, which is
based on the algorithm of
\cite{Majumder:1993:AAS} and the succesive improvements and corrections 
\cite{Cran:1977:RAS,Berry:1990:AAS,Berry:1991:CTA}. In \cite{Majumder:1993:AAS}, a simple approximation is used in
terms of the error function together with two additional starting value approximations for the tails; these initial values
are refined by the Newton iteration. As discussed in \cite{Cran:1977:RAS,Berry:1990:AAS}, the Newton method needs some
modification to ensure convergence inside the interval $[0,1]$, and further tuning of the Newton method has been considered
in recent versions of this algorithm for {\bf R} (but some convergence problem still remain in the present version, as
we later discuss).

In this paper the methods for the computation of the inversion of the cumulative beta distribution are improved in several
directions. In the first place, we introduce the Schwarzian-Newton method (SNM) as alternative to Newton's method (NM). 
With respect to Newton's method the SNM has the advantage of having order of convergence four instead of two.
In addition, as explained in \cite{Segura:2016:SNM}, the SNM has good non-local properties for this type of functions and
it is possible to build an algorithm with certified convergence. 
In the second place, we analyze initial value approximations (much sharper than those given in \cite{Majumder:1993:AAS})
in terms of asymptotic approximations for large values of $p$ and/or $q$, but which also give accurate values for
moderate values; these approximations are given in terms of inverse error functions or the inverse gamma distribution
(\cite{Temme:1992:AIB}, \cite[\S10.5.2]{Gil:2007:NSF}, 
\cite[\S42.3]{Temme:2014:AMI}). We also provide improved approximations for the tails obtained from the sharp bounds
described in \cite{Segura:2016:SBF}.

An additional direction of improvement of the algorithms is in the selection of the methods of
computation of the beta distribution, which are needed in the application of iterative methods (with Newton, SNM or any other
choice). This is not discussed in this paper, and we leave this topic for future research. A relatively recent algorithm was
given in \cite{Didonato:1992:A7S}.

\section{Methods of computation}

We next describe the methods of computation used in the algorithms. First we describe the SNM method,
and discuss how a standalone algorithm with certified convergence can be built with this method, provided an accurate
method of computation of the beta distribution is available. In the second place we describe the methods for estimating
the quantile function based on asymptotics for large $p$ and/or $q$. Finally, we describe sharp upper and lower bounds
for the tails that can be used for solving the problem (\ref{eq:int04}) for $\alpha$ close to zero or $1$.

\subsection{Schwarzian-Newton method}\label{sec:schw}

The Schwarzian-Newton method (SNM) is a fourth order fixed point
method with good non-local convergence 
properties for solving nonlinear equations $f(x)=0$ \cite{Segura:2016:SNM}. The SNM
 has Halley's method as limit when the Schwarzian derivative of the
function $f(x)$ tends to zero.

  Given a function $f(x)$ with positive derivative (in our case $f(x)=I_x (p,q)-\alpha$), 
it is easy to prove that $\Phi(x)=f(x)/\sqrt{f'(x)}$ satisfies the differential equation
\begin{equation}
\Phi'' (x)+ \Omega (x)\Phi (x)=0,\,\Omega=\Frac{1}{2}\{f,x\},
\end{equation}
where $\{f,x\}$ is the Schwarzian derivative of $f(x)$ with respect 
to $x$:

\begin{equation}
\{f,\,x\}= \Frac{f'''}{f'} -\Frac{3}{2}\left(\Frac{f''}{f'}\right)^2.
\end{equation}

  The SNM is obtained from the approximate integration of the Riccati equation $h'(x)=1+\Omega(x)h(x)^2$,
$h(x)=\Phi(x)/\Phi'(x)$ under the assumption that $\Omega(x)$ is approximately constant. In the
case of negative Schwarzian derivative (which will be the case for the beta distribution)
 the iteration function can be written as:

\begin{equation}
\label{atanh}
g(x)=x-\Frac{1}{\sqrt{|\Omega}|}{\rm arctanh}\left(\sqrt{|\Omega|} h(x)\right).
\end{equation}


 We discuss two implementations: a direct implementation, which gives a convergent algorithm for $p,q>1$ and
an implementation with an exponential change of variables, which is more easy to handle for the rest of cases.

\subsubsection{The direct implementation}

 It is proved in \cite{Segura:2016:SNM} that if $\Omega(x)$ has one and only one extremum at $x_e \in I$ 
and it is a maximum, then
 if $\Omega <0$ the SNM converges monotonically to the root of $f(x)$ in $I$ starting from $x_0=x_e$.
We will use this result for the cumulative central beta distribution, when the parameters $p$ and $q$ are
larger than 1. In this case, the function $\Omega (x)$ (the Schwarzian derivative of $f(x)$ with a factor $1/2$) is given by 

\begin{equation}
\label{omegabeta}
\Omega (x)=\Frac{(p-1)(q-1)}{2x(1-x)}-\Frac{1}{4}\Frac{p^2-1}{x^2}-
\Frac{1}{4}\Frac{q^2-1}{(1-x)^2}.
\end{equation}

It is possible to show that for $p>1$ and $q>1$, the function $\Omega(x)$ in (\ref{omegabeta})
is negative in $(0,1)$ and has only one extremum (which is a maximum) in that interval. The extremum
of $\Omega(x)$ is at

\begin{equation}
x_e=\Frac{1}{3\Delta^{1/3}}\Frac{(3pq+3p^2+6p)\Delta^{1/3}-\Delta^{2/3}+3pq\left((p+q)^2+8(p+q)+12\right)}{(p+q)^2+2(p+q)},
\end{equation}
where 

$$
\begin{array}{r@{\,}c@{\,}l}
\Delta&=&pq \left\{ 108(p-q)(p+q+1)+27(p^2q+p^3-q^2p-q^3)+3\sqrt{3}(p+q)\right.\\
       &&\left.(p+q+2)\displaystyle\sqrt{(p+q+2)(27q+54+q^2p+18pq+27p+p^2q)}\right\}.
\end{array}
$$

Then, the fixed point method is (\ref{atanh}) with $\Omega(x)$ given by (\ref{omegabeta}) and
$$
h(x)=\Frac{f(x)}{\Frac{1}{2}\left(-\Frac{p-1}{x}+\Frac{q-1}{1-x}\right)f(x)+
\Frac{x^{p-1}(1-x)^{q-1}}{B(p,q)}},\,
$$
where $f(x)=I_x (p,q)-\alpha$.

\subsubsection{The exponential implementation}

When $p$ and/or $q$ are smaller than 1, it is possible to make a change of variables in order to obtain
a negative Schwarzian derivative and simpler monotonicity properties. In particular, with the change of variables 

\begin{equation}
\label{changev}
z(x)=\log\left(\Frac{x}{1-x}\right),
\end{equation}
we obtain that $\Phi (z)=f(z)/\sqrt{\dot{f}(z)}$ (where the dot represents the derivative with respect to $z$) 
satisfies $\ddot{\Phi}(z)+\Omega(z)\Phi (z)=0$, where 
\begin{equation}
\label{omegabeta2}
\Omega (z)=\Frac{1}{4} \left(-(p+q)(p+q-2)x^2(z)+2(p+q)(p-1)x(z)-p^2 \right).
\end{equation}

The function $\Omega (z)$ has its extremum at $z_e=\log(x_e/(1-x_e))$, $x_e=(p-1)/(p+q-2)$.
When $p$ and/or $q$ smaller than 1, $\Omega(z)$ can be either be monotonic or it can have a minimum. 
Convergence of the SNM can be guaranteed, in this case, using the following
results \cite{Segura:2016:SNM}: a) if $\Omega(z)$ is negative and decreasing in the interval $I=[a,\,\alpha]$, then the SNM
converges monotonically to $\alpha$ for any starting value $z_0 \in [a,\,\alpha]$; b) if $\Omega(z)$ is negative
 and increasing in the interval $I=[\alpha,\,b]$, then the SNM
converges monotonically to $\alpha$ for any starting value $z_0 \in [\alpha,\,b]$.

The case $p=q=1$ is of course trivial. Apart from this, there are three different cases to be considered:

\begin{description}
\item{a)} $p\le 1,\,q>1$: the function $\Omega(z)$ is decreasing. In this
case, the SNM uses as starting point a large negative value (in the $z$ variable).  

\item{b)} $p>1,\,q\le 1$:  the function $\Omega(z)$ is increasing. In this
case, the SNM uses as starting point a large positive value (in the $z$ variable). 

\item{c)} $p<1,\,q<1$: the extremum of $\Omega(z)$ at $z_e$ is reached and it is a minimum. In this case,
we use the sign of the function $h(z)$ at $z_e$ to select a subinterval for application of the SNM, according
to the previous results. The function $h(z)$ is 
given by

\begin{equation}
\label{phixz}
h(z)=\Frac{f(x(z))}{-\Frac{1}{2}\Frac{p-qe^z}{1+e^z}f(x(z))+\Frac{e^{zp}}{B(p,q)(1+e^z)^{p+q}  }    } .
\end{equation}

When this sign is negative, the SNM uses a large positive value (in the $z$ variable), otherwise the SNM
uses a large negative value.

\end{description}

Once the SNM is applied to find the root $z_r$ in the $z$-variable, the corresponding $x$-value will be given by
$x_r=\Frac{e^{z_r}}{1+e^{z_r}}$.

\subsubsection{Discussion}

We have constructed two methods of order four which are proven to converge with certainty
for the initial values prescribed. The method has, in addition, good non-local properties, which means that
few iterations are needed for a good estimation of the inverse (typically from 2 to 4 for 20 digits), 
even without accurate starting values. The exceptions are the tails ($\alpha$ very close to $0$ or $1$), but
we will discuss later how to deal with these cases.

Because the convergence is guaranteed, no special precaution is needed to ensure that the interval $[0,1]$ 
in the original variable is not abandoned, as happened with earlier versions of the algorithm given in 
\cite{Majumder:1993:AAS} (see \cite{Cran:1977:RAS}) and as it is still happens for some values in the 
latests {\bf R} version of this algorithm. For instance, the {\bf R} command {\rm qbeta(alpha,600,1.1)} does
not converge properly if alpha$\in (6.9\,10^{-35},1.4\,10^{-20})$. Our method avoids this type of problems.

The performance of the method can be improved by considering initial approximations, which we are discussing next.

\subsection{Asymptotic inversion methods}\label{sec:asymp}

The algorithm considered in \cite{Majumder:1993:AAS}, which is the basis
of the {\bf R} implementation, uses an approximation in terms of the inverse error function, which 
works reasonably away from the tails. However, this simple approximation does not give more than 
two accurate digits, except by accident. 

Much more accurate initial approximations (some of them
also in terms of error functions) can be obtained from asymptotics for large $p$ and/or $q$.
These are accurate approximations for large and not so large $p$ and/or $q$, as we later discuss.

This section is based on the results given in \cite{Temme:1992:AIB} and 
\cite[\S42.3]{Temme:2014:AMI}.

\subsubsection{Inversion using the error function}\label{sec:asymperf}

We start with the following representation 

\begin{equation}\label{eq:pqerfc01}
I_x(p,q)=\tfrac12\erfc\left(-\eta\sqrt{r/2}\right)-R_r(\eta),
\end{equation}
where we write $p=r\sin^2 \theta$, $q=r\cos^2 \theta$ with $0<\theta <\pi/2$ and $\eta$ is given by

\begin{equation}\label{eq:pqerfc02}
-\tfrac12 \eta^2 = \sin^2 \theta \log \Frac{x}{\sin^2 \theta}+\cos^2 \theta \log \Frac{1-x}{\cos^2 \theta}. 
\end{equation}

When we take the square root for $\eta$, we choose $\sign(\eta)=\sign(x-\sin^2\theta)$, this means $\sign(\eta)=\sign(x-p/(p+q))$. In this way, the relation between $x\in(0,1)$ and $\eta\in(-\infty,\infty)$ becomes one-to-one.

Using this representation of $I_x(p,q)$, we solve the equation in \eqref{eq:int04} first in terms of $\eta$. When $r=p+q$ is a large parameter, the asymptotic method will provide an approximation to the requested value of $\eta$ in the form

\begin{equation}\label{eq:asinv03}
\eta \sim \eta_0 + \Frac{\eta_1}{r}+ \Frac{\eta_2}{r^2}+ \Frac{\eta_3}{r^3}+\ldots\,.
\end{equation}

The algorithm for computing the coefficients $\eta_i,\,\,i=0,1,2,\ldots$ can be summarized as follows

\begin{enumerate}
\item The value $\eta_0$ is obtained from the equation

\begin{equation}\label{eq:asinv04}
\tfrac12\erfc\left(-\eta_0\sqrt{r/2}\right)=\alpha.
\end{equation}

\item With $\eta=\eta_0$, equation (\ref{eq:pqerfc02}) is inverted to obtain a first approximation
to the value of $x$.
For inverting this equation, it seems convenient to write it in the form

\begin{equation}\label{eq:asinv05}
x^p(1-x)^q=\left(\Frac{p}{r} \right)^p \left(\Frac{q}{r} \right)^q e^{-r \eta^2/2}. 
\end{equation}

\item With these values of $\eta_0$ and $x$, the coefficient $\eta_1$ is given by

\begin{equation}\label{eq:asinv06}
\eta_1=\Frac{\log(f(\eta_0))}{\eta_0},
\end{equation}
where $f(\eta)=\Frac{\eta \sin\theta \cos \theta}{(x-\sin^2 \theta)}$.

\item Higher-order coefficients $\eta_j, \,j=2,3,\ldots$ can be obtained in terms of $x$, $\eta_0$, $\eta_1$, 
$\sin \theta$ and $\cos \theta$. As an example, the coefficient $\eta_2$ is given by
\begin{equation}\label{eq:asinv07}
\begin{array}{r@{\,}c@{\,}l}
\eta_2&=& \Frac{1}{12\eta_0^3c^2s^2(s^2-x)^2}
    \left(s^6\eta_0^2-\eta_0^2x^2-s^4\eta_0^2-\eta_0^2s^8\ +\right.\\[8pt]
   &&12s^6c^2-12s^2c^2\eta_1\eta_0^3x +12s^2c^2\eta_1\eta_0^3x^2-6\eta_0^2s^2c^2\eta_1^2x^2\ + \\[8pt]
   &&12\eta_0^2s^4c^2\eta_1^2x+2\eta_0^2xs^2+2\eta_0^2xs^6-6\eta_0^2s^6c^2\eta_1^2\ +\\[8pt]
   &&12s^2c^2\eta_0^2x^2-12s^2c^2\eta_0^2x-2\eta_0^2xs^4-\eta_0^2x^2s^4\ +\\[8pt]
    &&\left.\eta_0^2x^2s^2-24s^4c^2x+12s^2c^2x^2\right),
\end{array}
\end{equation}
where $s=\sin \theta$, $c=\cos \theta$.

\item
With these coefficients in the expansion (\ref{eq:asinv03}), a value for $\eta$ is obtained.
Then, the inversion of (\ref{eq:pqerfc02}) will provide the final asymptotic estimation of the $x$-value.  

\end{enumerate}

Using  \eqref{eq:pqerfc02} we can derive the following expansion or small values of $|\eta|$:
\begin{equation}\label{eq:asinv08}
\begin{array}{r@{\,}c@{\,}l}
x&=&\dsp{s^2+sc\,\eta+\frac{1-2s^2}3\eta^2+\frac{13s^4-13s^2+1}{36sc}\eta^3\ +}\\[8pt]
&&\dsp{\frac{46s^6-69s^4+21s^2+1}{270s^2c^2}\eta^4+\ldots,}
\end{array}
\end{equation}
where $s=\sin\theta$, $c=\cos\theta$. For larger values of $|\eta|$, with $\eta<0$, we rewrite \eqref{eq:pqerfc02} in
the form $x(1-x)^\mu=u$, where
\begin{equation}\label{eq:asinv09}
\mu=\cot^2\theta, \quad u=\exp\left[(-\tfrac12\eta^2+s^2\ln s^2+c^2\ln c^2)/s^2\right],
\end{equation}
and for small values of $u$ we expand
\begin{equation}\label{eq:asinv10}
\begin{array}{r@{\,}c@{\,}l}
x&=&\dsp{u +\mu u^2+\frac{3\mu (3\mu+1)}{3!}u^3+\frac{4\mu(4\mu+1)(4\mu+2)}{4!}u^4\ +}\\[8pt]
&&\dsp{\frac{5\mu (5\mu + 1)  (5\mu + 2)  (5\mu + 3)}{5!}u^5+ \ldots.}
\end{array}
\end{equation}

A similar approach is possible for positive values of $\eta$, giving an expansion for $x$ near
unity. In that case we have the equation $x^\nu(1-x)=v$, where
\begin{equation}\label{eq:asinv11}
\nu=\tan^2\theta, \quad v=\exp\left[(-\tfrac12\eta^2+s^2\ln s^2+c^2\ln c^2)/c^2\right],
\end{equation}
and we have the expansion
\begin{equation}\label{eq:asinv12}
\begin{array}{r@{\,}c@{\,}l}
1-x&=&\dsp{v +\nu v^2+\frac{3\nu (3\nu+1)}{3!}v^3+\frac{4\nu(4\nu+1)(4\nu+2)}{4!}v^4\ +}\\[8pt]
&&\dsp{\frac{5\nu (5\nu + 1)  (5\nu + 2)  (5\nu + 3)}{5!}v^5+ \ldots,}
\end{array}
\end{equation}

The approximations to $x$ obtained in this way will be used for starting the SNM for
obtaining  more accurate values of $x$.

\subsubsection{Inversion using the incomplete gamma function}

In this case, we start from

\begin{equation}\label{eq:invgam1}
I_x(p,q)=Q(q,\eta p)+R_{p,q}(\eta),
\end{equation}
where $Q(a,x)$ is the incomplete gamma function ratio.

The parameter $\eta$ is given by

\begin{equation}\label{eq:invgam2}
\eta -\mu\log\eta+(1+\mu)\log(1+\mu)-\mu=-\log x-\mu \log(1-x), 
\end{equation}
where $\mu=q/p$ and $x$ have the following corresponding points:
\begin{equation}\label{eq:invgam2a}
x= 0 \Longleftrightarrow \eta=+\infty\quad
x= 1/(1+\mu) \Longleftrightarrow \eta=\mu,\quad
x= 1 \Longleftrightarrow \eta=0.
\end{equation}
So, for $x\in(0,1)$ we have $\sign(\eta-\mu)=-\sign(x-1/(1+\mu))$.

 For the representation in \eqref{eq:invgam1} we assume that $p$ is the large parameter, and we will obtain approximations to the value
of $\eta$ in the form

\begin{equation}\label{eq:invgam3}
\eta \sim \eta_0 + \Frac{\eta_1}{p}+ \Frac{\eta_2}{p^2}+ \Frac{\eta_3}{p^3}+\ldots,\,\,p \rightarrow \infty.
\end{equation}
 
We follow similar ideas as in \S\ref{sec:asymperf}. The value of $\eta_0$ can be obtained by solving

\begin{equation}\label{eq:invgam3a}
Q(q,\eta_0 p)=\alpha.
\end{equation}

The inversion of $Q(a,x)$ can be done by using our inversion algorithm
described in \cite{Gil:2012:IGR}. 

Then, a value $x_0$ is obtained by solving (\ref{eq:invgam2}) for $x$. With $x_0$ and $\eta_0$
we compute

\begin{equation}\label{eq:invgam4}
\eta_1=\Frac{\log \phi(\eta_0)}{1-\mu/\eta_0},
\end{equation}
where $\phi(\eta)$ is given by $\phi(\eta)=\Frac{\eta-\mu}{1-x(1+\mu)}\Frac{1}{\sqrt{1+\mu}}$.
 
  Other coefficients $\eta_j,\,j=2,3,\ldots$ can be obtained in terms of $\mu$, $x_0$, $\eta_0$ and $\eta_1$.

To compute $x$ from equation \eqref{eq:invgam2} we can use the  inversion algorithm for computing $x$ when in \eqref{eq:pqerfc02} $\eta$ is given. This follows from $\mu=q/p=\cot^2\theta$ and from writing  \eqref{eq:invgam2} in the form
\begin{equation}\label{eq:invgam5}
\sin^2 \theta\,\left(\mu-\eta+\mu\log\frac{\eta}{\mu}\right)=
 \sin^2 \theta \log \frac{x}{\sin^2 \theta}+\cos^2 \theta \log \frac{1-x}{\cos^2 \theta}. 
\end{equation}
This equation can also be written as

\begin{equation}\label{eq:invgam6}
x^p(1-x)^q=\left(\Frac{p}{r} \right)^{r} e^{q(1+\log \eta) -p \eta}.   
\end{equation}

\subsection{Interval estimation for the tails}\label{sec:bounds}

 Sharp lower and upper bounds for the solution $x$ of the equation (\ref{eq:int04}) in the lower ($\alpha \rightarrow 0$)  and upper 
($\alpha \rightarrow 1$) tails 
of the distribution function 
can be obtained by using fixed point iterations $x^{n+1}_l=g_l(x^n_l)$ and $x^{n+1}_u=g_u(x^n_u)$, respectively, where 
the iteration functions $g_l(x)$ and $g_u(x)$ for the lower tail are given by \cite{Segura:2016:SBF}

\begin{equation}
g_l(x)=\left(\alpha B(p,q) \left(p-(p+q)x\right)(1-x)^{-q}\right)^{1/p},
\end{equation}
and

\begin{equation}
g_u(x)=\left(\Frac{\alpha p B(p,q)}{ \left(1+ \Frac{(p+q)}{(p+1)}x+  \Frac{(p+q)(p+q+1)}{(p+1)(p+2)}x^2\right)(1-x)^q}\right)^{1/p}.
\end{equation}

The starting value of the fixed point iterations is $x=0$. 
The solution $x$ of the equation (\ref{eq:int04}) satisfies $x_l<x<x_u$. These bounds of $x$ can be used 
as starting values for the SNM.
Notice that, because a lower and an upper bound is obtained we have an estimation of the error for these approximations that
can be used to decide if they are accurate enough.

We notice that the approximation for the lower tail used in \cite{Majumder:1993:AAS} (and also in the {\bf R} implementation) is
just $g_l (0)=g_u(0)$. Our approximation is more accurate and provides upper and lower bounds.

For the upper tail, the iteration functions are the same as before,  but with $p$ and $q$ interchanged  (by using (\ref{eq:int03})). 
The bounds are then given for $1-x$.

\section{Numerical testing}
 \label{sec:test}

In this section we illustrate the use of the different methods with numerical examples which
will help in deciding how to combine the methods in order to obtain fast reliable algorithms
(which are described in section \ref{algo}).

\subsection{Initial values obtained with the asymptotic approximations}
\label{subsec:testasymp}

In Tables 1 and 2 we show examples of the accuracy obtained in 
the computation of $\left|I_x(p,q)-\alpha \right|/\alpha$ with $x$ the values provided
by the asymptotic approximations (before iterating the SNM). \footnote{In Table 1 the
accuracy  $5.6\,10^{-16}$ corresponds to 
 the case $I_x(3,3)=0.5$: because of the
symmetry, $x$ should be $0.5$ (exact), and $\eta$
defined in (\ref{eq:pqerfc02}) becomes $0$, the same as $\eta_0$ in
(\ref{eq:asinv04}). This explains why that result in the table
becomes so small.}

The asymptotic methods provide relatively good initial values even
for quite small values of $p$ and $q$: using $10^7$ random points in the region 
$(p,q,\alpha) \in (0.5,\,1.5)\times (0.7,\,1.5)\times (0,\,1)$ we have tested that 
a relative accuracy better than $0.06$ was obtained when
computing $\left|I_x(p,q)-\alpha \right|/\alpha$ with $x$, the asymptotic
approximations obtained using the error function. With these initial values, not more than two iterations of the SNM are needed
to obtain an accuracy better than $5.0\,10^{-13}$. The function $I_x(p,q)$ is computed in the
iterations of the SNM by using a continued fraction representation 

\begin{equation}
\label{eq:contfrac}
I_x(p,q)=\Frac{x^p(1-x)^q}{pB(p,q)}
\left(\Frac{1}{1+}\,\Frac{d_1}{1+}\,\Frac{d_2}{1+}\,\Frac{d_3}{1+}\,\dots\right),
\end{equation}
where

\begin{equation}
\begin{array}{r@{\,}c@{\,}l}
d_{2m}&=& \Frac{m(q-m)x}{(p+2m-1)(p+2m)},\\
d_{2m+1}&=&-\Frac{(p+m)(p+q+m)x}{(p+2m)(p+2m+1)}.
\end{array}
\end{equation}

For the computation of the Beta function $B(p,q)$, it is convenient, in particular when $p$ and $q$ are large,
to use the following expression in terms of the scaled gamma function $\Gamma^*(x)$:

\begin{equation}
\label{eq:betaf}
B(p,q)=\sqrt{2 \pi} \displaystyle\sqrt{\Frac{1}{p}+\Frac{1}{q}}\left(\Frac{p^{\frac{p}{p+q}}q^{\frac{q}{p+q}}}{p+q}\right)^{p+q} 
\Frac{\Gamma^*(p) \Gamma^*(q)}
{\Gamma^*(p+q)},
\end{equation}
where $\Gamma^*(x)$ is defined as

\begin{equation}\label{eq:moc05}
\Gamma^*(x)=\frac{\Gamma(x)}{\sqrt{2\pi/x}\, x^xe^{-x}},\quad x>0.
\end{equation}

The function $\Gamma^*(x)$ is computed using the function {\bf gamstar}
included in a previous package developed by the authors \cite{Gil:2015:CPC}.

 \begin{table}
$$
\begin{array}{llll}
\hline
 \alpha & p=4 & p=3  & p=2   \\
  \hline
10^{-6} &  6.3\,10^{-4} & 1.6\,10^{-3} & 1.8\,10^{-3}  \\
10^{-3} &  3.2\,10^{-4} & 1.6\,10^{-3} & 4.5 \,10^{-3}  \\
0.1    & 2.7\,10^{-4}  &  4.0 \,10^{-4} & 1.9\,10^{-3}   \\ 
0.3    & 2.9 \,10^{-5}  & 3.9\,10^{-6} &  5.9\,10^{-5} \\
0.5    &  2.9\,10^{-5} & 5.6\,10^{-16} &  2.9 \,10^{-5}  \\
0.7    & 2.6\,10^{-5}  & 1.7\,10^{-6} &  1.2\,10^{-5}  \\
0.9    & 2.2\,10^{-4}  & 4.5\,10^{-5}  & 2.9\,10^{-5}  \\ 
0.999 &  4.5\,10^{-6}  &  1.6\,10^{-6} & 3.2\,10^{-7}  \\
0.99999 & 2.9\,10^{-8}  &  1.8\,10^{-8} &  6.2\,10^{-9}  \\
\end{array}
$$
{\footnotesize {\bf Table 1}. Relative errors 
$\left|I_x(p,q)-\alpha \right|/\alpha$ for $r=p+q=6$ using the estimates provided
by the asymptotic inversion method with the error function.   }
\label{table1}
\end{table}

 \begin{table}
$$
\begin{array}{lllll}
\hline
 \alpha & \mu=0.1 & \mu=0.5  & \mu=2  \\
  \hline
10^{-6} & 8.1\,10^{-5}  & 2.2\,10^{-4}   &   4.0 \,10^{-4}       \\
10^{-4} &   3.3\,10^{-4}  &  2.3\,10^{-5}  &  2.2\,10^{-4}         \\
0.1    &  2.4\,10^{-4}  &   1.9\,10^{-4}   &  2.5\,10^{-5}     \\ 
0.3    &  1.3\,10^{-4}   &  1.4\,10^{-4}   &   3.6\,10^{-5}       \\
0.5    &  7.8\,10^{-5}   &  9.8\,10^{-5}   &   3.2\,10^{-5}      \\
0.7    &  3.9\,10^{-5}   &   6.1\,10^{-5}  &  2.4\,10^{-5}        \\
0.9    &  1.1\,10^{-5}   &   2.3\,10^{-5}  &  1.1\,10^{-5}             \\
0.999 &   1.0\,10^{-7}   &      3.0\,10^{-7}   & 2.2\,10^{-7}   \\
0.99999 & 1.0\,10^{-9}    & 3.2\,10^{-9}     &2.9\,10^{-9}   \\
\end{array}
$$
{\footnotesize {\bf Table 2}. Relative errors 
$\left|I_x(p,q)-\alpha \right|/\alpha$ for $p=7$ and several values of $\mu=q/p$ using the estimates provided
by the asymptotic inversion method with the incomplete gamma function.   }
\label{table2}
\end{table}

\subsection{Initial values obtained in the tails of the distribution function}
\label{subsec:testbounds}

The interval estimations in the tails of the central beta
distribution function by using the fixed point iterations
 of \S\ref{sec:bounds}, can be also used  for providing starting values of the SNM.
This will be particularly useful for quite small values of the parameters $p$ and $q$, where the asymptotic method cannot be applied.
It is important to note that when the value of the parameters $p$ or $q$ are close to $0$, the inversion of $I_x(p,q)$ becomes
problematic in the lower (when $p \rightarrow 0$) or upper (when $q \rightarrow 0$) 
tail of the cumulative distribution function, because of the particular shape of the functions.

In Table 3 we show the relative errors $1-x_l/x$ and $1-x_u/x$ obtained with the lower and upper
 bounds, respectively, for the solution $x$ of the equation 
(\ref{eq:int04}) for small values of $p$, $q$ and $\alpha$.
The bounds (computed in the examples using Maple) 
are obtained  with just three iterations
of the fixed point methods  of \S\ref{sec:bounds}. We have also tested that for small values
of $p$ and $q$, the bounds provide in all cases reasonable approximations for starting the SNM, no matter
if the value of $\alpha$ is small. 
Besides, even for not so small values of $p$ and $q$, the bounds
provide very accurate estimations when $\alpha$ is very small. In some cases, these estimations could be even better
than the estimations of the asymptotic method.

 \begin{table}
$$
\begin{array}{llcc}
\hline
 \alpha &  &p=0.3   &  p=0.4     \\
        &  & q=0.4   &  q=0.3     \\
  \hline
10^{-7}  & \mbox{(LB)}  & -1.2\,10^{-22} &      -5.9\,10^{-17}  \\
         &  \mbox{(UB)} &  7.8\,10^{-49}   &     1.1\,10^{-49}  \\
\hline 
10^{-5}  & \mbox{(LB)}  & -5.4\,10^{-16}      &  -5.9\,10^{-12}   \\
         &  \mbox{(UB)} &   1.3\,10^{-48}     &   5.6\,10^{-36}    \\
\hline 
10^{-3}  & \mbox{(LB)}  &  -2.5\,10^{-9}      &   -5.9\,10^{-7}   \\
         &  \mbox{(UB)} &    8.5\,10^{-29}     &   5.6\,10^{-21}    \\
\hline
\end{array}
$$
{\footnotesize {\bf Table 3}. Relative errors 
 $1-x_l/x$ and $1-x_u/x$ obtained with the lower (LB) and upper (UB) bounds for the solution $x$ of the equation 
(\ref{eq:int04}). }
\label{table3}
\end{table}

\subsection{Performance of the SNM for small values of the parameters }

We have tested that the scheme for the SNM, as described in  \S\ref{sec:schw} when the parameters $p$ and $q$ are both small,
also provides a good uniform accuracy in the computation: 
using $10^7$ random points in the region $(p,q,\alpha) \in (0.1,\,0.5)\times (0.1,\,0.7) \times (0,\,1)$
we have tested that 
a relative accuracy better than $4.8\,10^{-13}$ was obtained when
computing $\left|I_x(p,q)-\alpha \right|/\alpha$. The maximum number of iterations of the SNM was 3.

\subsection{Efficiency testing}

As we have shown in \S\ref{subsec:testasymp}, the asymptotic method provides very accurate initial values for starting the
SNM even for small values of the
parameters $p$ and $q$. But apart from accuracy, an important feature of any computational scheme is also efficiency. 
So, we have compared whether the combined use of the asymptotic approximations plus iterations of the SNM
is more efficient or not than the sole use of the SNM. In Table 4 we show CPU times spent by 20000 runs
of the inversion algorithm for different values of $p$, $q$ and $\alpha$ using three methods
of computation:  a) asymptotic inversion 
method using the error function for estimating the initial value plus iterations of the SNM; b)
 asymptotic inversion method using the gamma function for estimating the initial value plus iteration of the
SNM; c) iterations of the SNM with starting values obtained as discussed in \S\ref{sec:schw}.  In all cases the SNM is iterated until 
the solution of equation (\ref{eq:int04}) is obtained with an accuracy near full double precision. 

 \begin{table}
$$
\begin{array}{llccccc}
\hline
 \alpha & \mbox{Method}&p=4  & p=50, & p=100,  &p=150, &  p=300, \\
        &              & q=3  &q=60 & q=80  &q=1.0 &  q=400  \\
\hline
10^{-6} & \mbox{M1} & 0.22   &0.29     &   0.22  & 0.4   & 0.20   \\
        &\mbox{M2}  &  0.31  &0.39    &    0.31 & 0.22 &   0.33      \\
        &\mbox{M3}  & 0.32  & 0.34    &    0.34  & 0.39 &   0.38        \\
\hline
10^{-4} & \mbox{M1} & 0.22  &0.23     &   0.19  & 0.36   &  0.22  \\
        &\mbox{M2}  & 0.34  &0.34    &    0.31 & 0.33 &    0.33     \\
        &\mbox{M3}  & 0.25  &0.27    &    0.28  & 0.27 &   0.33     \\
\hline
0.3     & \mbox{M1} & 0.19  & 0.16   &  0.17   &  0.20     &    0.17   \\
        & \mbox{M2} & 0.31  &0.28    &   0.30    &  0.19    &  0.30     \\
        & \mbox{M3} & 0.11  &0.19   &   0.20    &  0.19    &  0.17        \\
\hline
0.7    & \mbox{M1} & 0.20   &0.17   &  0.16   &  0.23    &   0.18    \\
        & \mbox{M2} & 0.33  &0.28  &   0.30   &  0.23   &  0.30      \\
        & \mbox{M3} & 0.16  & 0.19  & 0.19     &  0.16     &  0.19       \\
\hline
0.999   & \mbox{M1} & 0.17  &0.17     &  0.17   &   0.14   &  0.16     \\
        & \mbox{M2} & 0.25 &0.28      &  0.28    &  0.16    &  0.33       \\
        & \mbox{M3} & 0.25 &  0.27   &  0.28     & 0.16      &    0.27     \\
\end{array}
$$
{\footnotesize {\bf Table 4}. CPU times (in seconds) for 20000 runs
of the inversion algorithm using different methods. M1: Asymptotic inversion 
method using the error function +SNM; M2: Asymptotic inversion 
method using the gamma function +SNM; M3: SNM. 
 The SNM is iterated until the solution of equation (\ref{eq:int04}) is obtained with an accuracy 
near full double precision.  }
\label{table4}
\end{table}

The results in Table 4 and additional testing for other parameter values,
indicate that the sole use of the SNM is efficient in all cases for the inversion of the
cumulative central beta distribution, but specially when the values of the parameter
$\alpha$ are neither very small nor near to $1$. 

 A different scenario arises when the solution of equation (\ref{eq:int04}) is computed
with an aimed accuracy better than $10^{-8}$ (single precision). In this case, just using the approximations
provided by the asymptotic expansions will be enough to obtain such an accuracy for a wide range
of parameters.

\section{Proposed algorithms}
\label{algo}

Based on the previous numerical experiments of \S\ref{sec:test}, we conclude that if the precision required is not
very high, the initial approximations given by asymptotics or by the tail estimations could
be sufficient in a large range of the parameters. However, for higher precision the use
of the SNM must be prevalent. 

This leads us to suggest two different schemes for computing the solution $x$ of the equation (\ref{eq:int04}).

\begin{description}
\item{SCHEME 1.}  Algorithm for the inversion of the cumulative central beta distribution with an accuracy
near double precision:

\begin{description}

\item{If $\alpha \le 0.01$}

\begin{description} 
\item{For} $p<0.3$, use the upper bound of \S\ref{sec:bounds} as solution
of the equation (\ref{eq:int04}). 

\item{For} $0.3<p<1$, use the SNM as described in \S\ref{sec:schw}
using as starting values the bounds of \S\ref{sec:bounds}.
\item{For} $1<p<30$ and $q<1$, use the SNM, using as starting values the bounds of \S\ref{sec:bounds}.

 \item{For} $p>30$ and $q<0.5$, use the SNM, using as starting values the bounds of \S\ref{sec:bounds}.

\item{For} $p>30$, $0.5<q<5$: a) if $\alpha>0.0001$ use the SNM, using as starting values the  
approximation provided by the uniform asymptotic expansion 
       in terms of the gamma function; b) if $\alpha<0.0001$ use the SNM, using as starting values 
the bounds of \S\ref{sec:bounds}.

 \item{In other cases,} use the SNM as described, using as starting values 
the approximations provided by the uniform asymptotic expansion 
       in terms of the error function in other cases.
\end{description} 
 
\item{When $0.01<\alpha \le 0.5$}

\begin{description} 
\item{For} $1<q<5$ and $p>50$,  use the SNM, using as starting values the 
approximations provided by the uniform asymptotic expansion 
       in terms of the incomplete gamma function. 
\item{For} $p>30$ and $q>30$,  use the SNM, using as starting values 
the approximations provided by the uniform asymptotic expansion 
       in terms of the error function. 
\item{In other cases,} use the SNM as described in \S\ref{sec:schw}.

\end{description}

\item{For $0.5<\alpha <1$}, use the relation (\ref{eq:int03}) and apply the previous steps
to solve $1-x$ in $I_{1-x}(q,p)=1-\alpha$.

\end{description}

\item{SCHEME 2.}  Algorithm for the inversion of the cumulative central beta distribution with an accuracy
near single precision:

\begin{description}
\item{If $\alpha \le 0.01$}

\begin{description} 
\item{For} $p<0.5$, use the upper bound of \S\ref{sec:bounds} as solution
of the equation (\ref{eq:int04}). 

\item{For} $0.5<p<1$, use the SNM as described in \S\ref{sec:schw}
using as starting values the bounds of \S\ref{sec:bounds}.

 \item{For} $1<p<30$ and $q<1$, use the SNM, using as starting values the bounds of \S\ref{sec:bounds}.

 \item{For} $p>30$ and $q<0.5$, use the SNM, using as starting values the bounds of \S\ref{sec:bounds}.

\item{For} $p>30$, $0.5<q<5$: a) if $\alpha>0.0001$ use the SNM, using as starting values the  
approximation provided by the uniform asymptotic expansion 
       in terms of the gamma function; b) if $\alpha<0.0001$ use the SNM, using as starting values 
the bounds of \S\ref{sec:bounds}.

 \item{In other cases,} use the SNM as described in \S\ref{sec:schw}
using as starting values the approximations provided by the uniform asymptotic expansion 
       in terms of the error function in other cases.
\end{description} 

\item{When $0.01<\alpha \le 0.5$}

\begin{description} 

\item{For} $1<q<3$, $p>160$ and $\alpha>0.1$, use the approximation provided by the uniform asymptotic expansion 
       in terms of the incomplete gamma function as solution
of the equation (\ref{eq:int04}). 
\item{For} $p>30$ and $q>30$, use the approximation provided by the uniform asymptotic expansion 
       in terms of the error function as solution
of the equation (\ref{eq:int04}). 
\item{In other cases,} use the SNM as described in \S\ref{sec:schw}.
\end{description}

\item{For $0.5<\alpha <1$}, use the relation (\ref{eq:int03}) and apply the previous steps
to solve $1-x$ in $I_{1-x}(q,p)=1-\alpha$.

\end{description}

\end{description}

\section{Conclusions}

We have presented methods for the inversion of cumulative beta distributions which improve 
previously existing methods. We have described how the Schwarzian-Newton method provides
a standalone method with certified convergence which is in itself an efficient method, even without accurate initial
estimations (except at the tails). In addition, we have discussed how to improve the efficiency by estimating
the quantile function using asymptotics for large $p$ and/or $q$ and by considering sharp upper and lower bounds
for the tails. These initial estimations are considerably more accurate than the simple approximations used in
some standard mathematical software packages (like {\bf R}) and, combined with the fourth order SNM, provide
efficient and reliable algorithms for the inversion of cumulative beta distributions.

\section{Acknowledgements}

The authors acknowledge financial support from 
{\emph{Ministerio de Econom\'{\i}a y Competitividad}}, project MTM2012-34787. NMT thanks CWI, Amsterdam, for scientific support.

\bibliographystyle{plain}
\bibliography{biblio}

\end{document}